\documentclass[runningheads]{llncs}

\usepackage{graphicx} % Required for inserting images
\usepackage{amsmath}
\usepackage{amsfonts}
\usepackage{amssymb}
\usepackage[lined,boxed,commentsnumbered,ruled]{algorithm2e}
\usepackage{geometry}
\usepackage{comment}

\usepackage{xcolor}

\newcommand{\sign}[1]{\mathrm{sgn}(#1)}

\title{A study on two-metric projection methods}
\author{Hanju Wu \and Yue Xie}
\institute{The University of Hong Kong}
\date{}

\begin{document}

\maketitle

\section{Abstract}
The two-metric projection method is a simple yet elegant algorithm proposed by Bertsekas in 1984 to address bound/box-constrained optimization problems. 
The algorithm's low per-iteration cost and potential for using Hessian information makes it a favourable computation method for this problem class. 
However, its global convergence guarantee is not studied in the nonconvex regime. 
In our work, we first investigate the global complexity of such a method for finding first-order stationary solution. 
After properly scaling each step, we equip the algorithm with competitive complexity guarantees. 
Furthermore, we generalize the two-metric projection method for solving $\ell_1$-norm minimization and discuss its properties via theoretical statements and numerical experiments.

\section{Introduction}

In this work, we revisit the two-metric projection method proposed by Bertsekas \cite{bertsekas1982projected} to solve the bound-constrained problem:
\begin{equation}\label{opt: bdcons}
    \begin{array}{*{2}{rl}}
        \min \; f(x) \quad \mbox{subject to} \quad x \geq 0
    \end{array}
\end{equation}
where  $f: \mathbb{R}^{n} \rightarrow \mathbb{R}$ is Lipschitz continuously differentiable 
and is bounded below by  $f_{\text {low }}$  on the feasible region. Two-metric projection method is simple and elegant, but its global complexity guarantees are absent in literature. We observe that this traditional method finds an approximate first-order optimal point in $\mathcal{O}(\epsilon^{-3})$ iterations, which is suboptimal in literature \cite{carmon2020lower}. By properly scaling the diagonal matrix, we equip the method with a worst-case complexity that matches the lower bound ($\mathcal{O}(\epsilon^{-2})$) of using first-order methods to solve this problem class. 
%Bounds represent the most basic form of inequality constraints. Subproblems with bound constraints frequently emerge in algorithms designed for more complex constrained optimization problems, including augmented Lagrangian methods. Additionally, bound constraints are commonly encountered in machine learning tasks such as nonnegative least-squares and nonnegative matrix factorization \cite{wang2012nonnegative}.
Moreover, we try to generalize the method and propose a two-metric adaptive projection method for the $\ell_1$-norm minimization,
\begin{equation}\label{opt: l1}
        \min \; \psi(x) \triangleq f(x)+\gamma \|x\|_1.
\end{equation}
We show that our method is well-defined and decreases the value of the objective function at each iteration $k$ for which $x_k$ is not a critical point.
For \eqref{opt: l1}, inexact proximal-Newton-type methods, or successive quadratic approximation (SQA) methods \cite{nocedal1999numerical}, are widely used for their superlinear
convergence in terms of iterations. A key difference between the two-metric adaptive projection and these methods is the per-iteration cost: in each iteration, SQA needs to find approximate solutions to regularized quadratic programs, while our method at worst needs to solve a linear equation - the Newton system and compute a single cheap adaptive projection step. Preliminary numerical experiments verify the asymptotic superlinear convergence of the two-metric adaptive projection method when employing Hessian.
% However, they suffer from lengthy running times in solving regularized subproblems because even approximate solutions cannot be computed easily. Our method, however, only needs to calculate a "projection" step at each iteration which is simpler than these methods.
\section{First-order optimal point}
%\subsection{First-order optimal points}
%\begin{theorem}[First-order necessary conditions]
We can write first-order optimality conditions for (1), (also known as stationarity conditions) at a point $\bar{x}$ as follows:
\begin{equation}
        \begin{aligned}
            \partial f\left(\bar{x}\right) / \partial x^{i}=0 \quad  \text{if}  \quad \bar{x}^{i}>0, \quad i=1, \ldots, n ,\\
            \partial f\left(\bar{x}\right) / \partial x^{i} \geq 0 \quad  \text{if}  \quad \bar{x}^{i}=0, \quad i=1, \ldots, n .
        \end{aligned}
    \end{equation}
%\end{theorem}

\subsection{Approximate first-order optimal points}
Next, we define $\epsilon$-approximate first-order optimal points ($\epsilon$-1o). %Our definition of an $\epsilon$-1o point is as follows.

\begin{definition}\label{def: 1o}
    For $\epsilon>0$, $\bar{x}$ is an $\epsilon$-1o point of (1), if $\bar{x} \geq 0$ and
    \begin{equation}
        \|S \nabla f(\bar{x})\| \leq \epsilon, \quad(\nabla f(\bar{x}))_{i} \geq-\epsilon,
    \end{equation}
    where $S$ is a diagonal matrix such that $S[i, i]=\min \left\{\bar{x}^{i}, 1\right\}$  if  $\bar{x}^{i} \geq 0,(\nabla f(\bar{x}))_{i}>0$; $S[i, i]=1$ otherwise.
\end{definition}

Definition~\ref{def: 1o} is motivated by the first-order optimal conditions of (1). 
In fact, if we let $\epsilon$ = 0, then 0-1o satisfies (3) exactly.
The following lemma further justifies Definition~\ref{def: 1o} and our purpose to find an
$\epsilon$-1o given small $\epsilon$.
\begin{lemma}
    Consider problem (1), suppose that we have a positive scalar sequence $\left\{\epsilon_k\right\}$ with 
    $\epsilon_k \to 0$ and vector sequence $\left\{x_k\right\}$ with $x_k\geq0$ and $x_k \to x_*$ such that $x_k$ is $\epsilon_k$-1o according to Definition 1.
    Then $x_*$ satisfies first-order optimal conditions (3).
\end{lemma}

\begin{proof}
    Denote diagonal matrix $S_k$ where coorespond to $S$ in Definition 1 with $x=x_k,\epsilon=\epsilon_k$.
    Our claim that $x_*$ satisfies (3) is a consequence of the following three observations.

    \begin{enumerate}
        \item Feasibility of $x_*$ follows from closedness of $\Omega$.
        \item $\forall i,k$, $ \nabla_i f(x_k) \geq -\epsilon_k$. By taking limits, we have $\nabla_i f(x_*) \geq 0$
        \item We claim that if $x_*^i > 0$, we must have $\nabla_i f(x_*) = 0$. Assume the contrary, if $\nabla_i f(x_*) > 0$, there exist a $K$, for all $k \geq K$ we have $\nabla_i f(x_k)>0$ and $x_k^i > \frac{1}{2}x_*^i > 0$, then $\forall k \geq K$
        \begin{equation*}
           \min \left\{\frac{1}{2}x_*^i, 1\right\}\nabla_i f(x_k) < |\min \left\{x_k^{i}, 1\right\}\nabla_i f(x_k)| \leq \|S_k \nabla f(x_k)\| \leq \epsilon_k.
        \end{equation*}
        By taking limits, we have $\nabla_i f(x_*)=0$, which is contradicting with $\nabla_i f(x_*) > 0$.
        % Consider $\nabla_i f(x_*) < 0$, there exist a $K$, for all $k \geq K$ we have $\nabla_i f(x_k)<0$, then
        % \begin{equation*}
        %     |\nabla_i f(x_k)| \leq \|S_k \nabla f(x_k)\| \leq \epsilon_k.
        % \end{equation*}
        % By taking limits, we have $\nabla_i f(x_*)=0$, which is contradicting with $\nabla_i f(x_*) < 0$.
    \end{enumerate}
    Overall, $x^*$ satisfies (3).
\end{proof}

\section{Two-metric projection}
The two-metric projection method for solving \eqref{opt: bdcons} is as follows:
$$x_{k+1}:=\mathcal{P}\left(x_{k}-\alpha_{k} D_{k} g_{k}\right), \quad \forall k \geq 0,$$
where $\mathcal{P}$ denotes Euclidean projection onto the feasible region, $\alpha_{k}>0$ is the stepsize, 
$g_{k} \triangleq \nabla f\left(x_{k}\right)$. 
$D_{k}$ is a symmetric positive definite matrix in $\mathbb{R}^{n \times n}$.
Define 
$$
p_{k} \triangleq D_{k} g_{k}, I_{k}^{+} \triangleq\left\{i\mid 0 \leq x_{k}^{i}\leq \epsilon_k, \;g_{k}^{i}>0\right\}, \; I_{k}^{-} \triangleq[1, \ldots, n] \backslash I_{k}^{+},
$$
where $\epsilon_k=\min \left\{\epsilon,\omega_k\right\}$, $\omega_k=\|x_k-P(x_k-M \nabla f(x_k))\|$, $\epsilon>0$ is a fixed tolerance and $M$ is a fixed diagonal positive definite matrix.
We suppose that  $D_{k}$  is diagonal w.r.t.  $I_{k}^{+}$ , i.e.,
\begin{equation}
    D_{k}[i, j]=0, \quad \forall i \in I_{k}^{+}, j \neq i .
\end{equation}
Denote $x_{k}(\alpha) \triangleq P\left(x_{k}-\alpha p_{k}\right), \forall \alpha>0$.
Let  $\bar{g}_{k}$  and  $\bar{p}_{k}$  be the subvector of  $g_{k}$  and  $p_{k}$  w.r.t index set  $I_{k}^{-}$ , and  $\bar{D}_{k}$ be the submatrix of $D_{k}$ and w.r.t.  $I_{k}^{-} \times I_{k}^{-}$. 
$S_k$ is the diagonal matrix which corresponds to $I^+, I^-$, and $S$ in Definition 1 with $x=x_k, s=s_k$ and $\epsilon$.
Then the algorithm can be formally stated as follows.
\begin{algorithm}\label{alg: 2mproj}
    \SetAlgoLined
    \KwIn{tolerance $\epsilon>0$, step acceptance parameter $\sigma \in(0,1)$, backtracking parameter $\beta \in(0,1)$ and initial feasible point $x_{0}$. Let $k=0$.}%输入参数
    \While{$\|S_k \nabla f(x_k)\|>\epsilon$}
    {
        Find:
        \begin{equation*}
            m_k \triangleq \min \left\{m \in \mathbb{N} \Bigg| f\left(x_{k}\right)-f\left(x_{k}\left(\beta^{m}\right)\right) \geq \sigma\left\{\beta^{m} \sum_{i \notin I_{k}^{+}} g_k^i p_{k}^{i}
        +\sum_{i \in I_{k}^+} g_k^i (x_{k}^{i}-x_{k}^{i}\left(\beta^{m}\right))\right\}\right\}
        \end{equation*}
        
        Let $\alpha_k= \beta^{m_k}$, $x_{k+1}= x_k(\alpha_k)$.
        
        Let $k=k+1$.
    }
    \KwOut{$x_k$}
    \caption{Two-metric projection}
\end{algorithm}
\begin{remark}
    We do not need to consider if $( \nabla f(x_k))_{i} < -\epsilon$ or not in the Algorithm 1, since if $(\nabla f(x_k))_{i} < -\epsilon$, then
    \begin{equation*}
        \|S_k \nabla f(x_k)\| \geq | (\nabla f(\bar{x}))_{i}| > \epsilon
    \end{equation*}
\end{remark}
Complexity of Algorithm~\ref{alg: 2mproj} is revealed in the following theorem.
\begin{theorem}\label{thm: 2mcomplexity}
    Consider Algorithm 1. Suppose that $\lambda_{\min }\|z\|^{2} \leq z^{T} D_{k} z \leq \lambda_{\max }\|z\|^{2}$, for all $z \in \mathbb{R}^{n}$ and $k \geq 0$. $\left\|g_{k}\right\| \leq G, \forall k \geq 0$. Lipschitz constant for  $\nabla f$  is $L$. 
    Then  $m_{k}<+\infty$  and the following holds: for any $0<\epsilon<1$, the algorithm will stop within $\mathcal{O}(\epsilon^{-3})$
    number of iterations and output an $\epsilon$-1o $\bar{x}$ that satisfies Definition 1.
\end{theorem}
\begin{remark}
We suppress the proof of Theorem~\ref{thm: 2mcomplexity} due to page limit. Note that the best complexity of first-order methods for Lipschitz smooth function $f$ in literature is $\mathcal{O}(\epsilon^{-2})$. The reason to cause this gap is the possible small stepsize in the current form of two-metric projection. Motivated by this gap, we design a scaled two-metric projection method to enlarge the stepsize and improve the complexity guarantees. 
%{\yx Does small stepsize really happen in reality? Is it only a problem in theory?} 
\end{remark}

\section{Scaled two-metric projection}
In this section, we apply an extended scaled version of the
two-metric projection algorithm. The only difference between the scaled version and the traditional version is that we use a diagonal matrix to scale $D_k$.
We update $x_{k+1}$ using:
$$x_{k+1}:=\mathcal{P}\left(x_{k}-\alpha_{k} S_{k} D_{k} S_{k} g_{k}\right), \quad \forall k \geq 0,$$
Where $S_{k}$ is defined as follows:
\begin{equation}
    S_{k}[i, i]=\min \left\{x_{k}^{i}, 1\right\} , \text{ if } g_k^i>0 ; \quad S_{k}[i, i]=1 \text {, otherwise. }
\end{equation}
Define $p_{k} \triangleq S_{k} D_{k} S_{k} g_{k}$. The Algorithm is formally the same as Algorithm 1 with only different $p_k$. The following theorem shows an improved complexity guarantee for this algorithm. 
\begin{theorem}
    Consider Algorithm 2. Suppose that $\lambda_{\min }\|z\|^{2} \leq z^{T} D_{k} z \leq \lambda_{\max }\|z\|^{2}$, for all $z \in \mathbb{R}^{n}$ and $k \geq 0$. $\left\|g_{k}\right\| \leq G, \forall k \geq 0$. Lipschitz constant for  $\nabla f$  is $L$. 
    Then  $m_{k}<+\infty$  and the following holds. 
    For any $0<\epsilon<1$, the algorithm will stop within $\left\lceil\frac{(f(x_0)-f_{\text{low}})\max \left\{4\lambda_{\max } G, \frac{2L \lambda_{\max }}{1-\sigma},2\beta \lambda_{\text{min}}\right\}}{\sigma \beta \lambda_{\text{min}} \epsilon^2}\right\rceil$ 
    number of iterations and output an $\bar{x}$ that satisfies Definition 1.
\end{theorem}

\begin{proof}
    First we show that  $m_{k}<+\infty, \forall k \geq 0$. 
    Denote $I_{k}^{p+} \triangleq\left\{i \mid p_{k}^{i}>0, g_{k}^{i} \leq 0\right\}$, $s_{k}^{i} \triangleq S_{k}[i, i]$. 
    Note that $I_{k}^{p+} \in I_k^-$ and $I_k^- \setminus I_{k}^{p+}=\left\{i \mid 0 \leq x_k^i \leq \epsilon_k, g_k^i \leq 0, p_k^i \leq 0\right\} \cup \left\{i \mid x_k^i > \epsilon_k, g_k^i > 0\right\} \cup \left\{i \mid x_k^i > \epsilon_k, p_k^i \leq 0\right\}$.
    If $\forall \alpha>0$ such that
    \begin{equation}
        \alpha \leq \min _{i \in\left\{i \in I_k^-\mid g_{k}^{i}>0, p_{k}^{i}>0\right\}}
        \left\{\frac{x_{k}^{i}}{p_{k}^{i}}\right\},
    \end{equation}
    we have the following properties:
    \begin{equation}
        x_k^i(\alpha)=x_k^i-\alpha p_k^i, \quad
        \forall i \in I_k^- \setminus I_{k}^{p+},
    \end{equation}
    \begin{equation}
        g_k^i \leq 0, x_k^i(\alpha)=\left[x_k^i-\alpha p_k^i\right]_+ \geq x_k^i-\alpha p_k^i \Longrightarrow g_k^i(x_k^i-\alpha p_k^i-x_k^i(\alpha)) \geq 0
        ,\quad \forall i \in I_{k}^{p+},
    \end{equation}
    Therefore,
    \begin{equation}
        \begin{aligned}
            & \sum_{i \in I_k^-} g_{k}^{i}\left(x_{k}^i(\alpha)-x_{k}^i\right) \\ 
            = & \sum_{i \in I_{k}^{-} \setminus I_k^{p+}} g_{k}^{i}\left(x_{k}^{i}(\alpha)-x_{k}^{i}\right)+\sum_{i \in I_{k}^{p+}} g_{k}^{i}\left(x_{k}^{i}(\alpha)-x_{k}^{i}\right) \\
            = &\sum_{i \in I_{k}^{-}\setminus I_k^{p+}} g_{k}^{i}\left(-\alpha p_{k}^{i}\right)+\sum_{i \in I_{k}^{p+}} g_{k}^{i}\left(x_{k}^{i}(\alpha)-x_{k}^{i}+\alpha p_{k}^{i}\right)+\sum_{i \in I_{k}^{p+}} g_{k}^{i}\left(-\alpha p_k^i\right) \leq  \sum_{i \in I_{k}^{-}}g_{k}^{i}\left(-\alpha p_{k}^{i}\right)
        \end{aligned}
    \end{equation}
    Since $\forall \alpha >0$, we have:
    \begin{equation}
        x_k^i(\alpha)=\left[x_k^i-\alpha p_k^i\right]_+ \geq x_k^i-\alpha p_k^i \Longrightarrow x_k^i-x_k^i(\alpha) \leq \alpha p_k^i, \quad \forall i \in I_k^+,
    \end{equation}
    \begin{equation}
        0 \leq p_k^i= s_k^i D_k\left[i,i\right] s_k^i g_k^i \leq D_k\left[i,i\right] g_k^i \leq \lambda_{\text{max}} g_k^i \Longrightarrow 0\leq x_k^i-x_k^i(\alpha), \quad \forall i \in I_k^+,
    \end{equation}
    \begin{eqnarray}
        \left | x_k^i-x_k^i(\alpha) \right | \leq \alpha \left | p_k^i \right | .
    \end{eqnarray}
    Then,
    \begin{equation}
        \sum_{i \in I_k^+} \left(x_{k}^i(\alpha)-x_{k}^i\right)^2
        \leq \sum_{i \in I_{k}^{+}} \alpha p_{k}^{i} (x_{k}^i-x_{k}^i(\alpha))
        \leq \alpha \lambda_{\text{max}} \sum_{i \in I_{k}^{+}}  g_{k}^{i} (x_{k}^i-x_{k}^i(\alpha)).
    \end{equation}
    Note that $\| \bar{S}_{k} \| \leq  1$ and $\forall z \in \mathbb{R}^n, D \in \mathbb{S}^{n}_{++}$
    we have:
    \begin{eqnarray}
        \|Dz\|^2 = z^{T}D^2 z = z^{T}D^{\frac{1}{2}} D D^{\frac{1}{2}} z = (D^{\frac{1}{2}} z)^{T} D (D^{\frac{1}{2}} z) \leq \lambda_{\text{max}} \|D^{\frac{1}{2}} z\|^2 = \lambda_{\text{max}} z^{T} D z
    \end{eqnarray}
    Therefore,
    \begin{equation}
        \begin{aligned}
            & \sum_{i \in I_k^-} \left(x_{k}^i(\alpha)-x_{k}^i\right)^2
            \leq \sum_{i \in I_{k}^{-}} (\alpha p_{k}^{i})^2
            = \alpha^2 \|\bar{p}_{k}\|^2 
            = \alpha^{2} \left\|\bar{S}_k \bar{D}_{k} \bar{S}_k \bar{g}_{k}\right\|^{2} \\
            \leq & \alpha^{2} \| \bar{S}_k \|^2 \left\|\bar{D}_{k} \bar{S}_k \bar{g}_{k}\right\|^{2} 
            \leq \alpha^{2}\left\|\bar{D}_{k} \bar{S}_k \bar{g}_{k}\right\|^{2}
            \leq \alpha^{2} \lambda_{\max } \bar{g}_{k}^{T} \bar{S}_k \bar{D}_{k} \bar{S}_k \bar{g}_{k} \\
            = & \alpha^{2} \lambda_{\max } \bar{g}_{k}^{T} \bar{p}_{k} =  \alpha^2 \lambda_{\text{max}} \sum_{i \in I_k^-} g_k^i p_k^i.
        \end{aligned}
    \end{equation}
    Since $\forall i \in \left\{i \in I_k^-\mid g_{k}^{i}>0, p_{k}^{i}>0\right\}$, we have:
    \begin{equation}
        \begin{aligned}
            0 < p_k^i = s_k^i (D_k S_k g_k)^i =\min \left\{x_k^i,1\right\}(D_k S_k g_k)^i \leq x_k^i \|D_k S_k g_k\|_{\infty}\\
            \leq x_k^i \|D_k S_k g_k\|_{2} \leq x_k^i \|D_k\|_2 \|S_k\|_2 \|g_k\|_2 \leq x_k^i \lambda_{\text{max}} G
        \end{aligned}
    \end{equation}
    which indicates that,
    \begin{equation}
        \frac{x_k^i}{p_k^i} \geq \frac{1}{\lambda_{\text{max}} G}, \quad \forall i \in \left\{i \in I_k^-\mid g_{k}^{i}>0, p_{k}^{i}>0\right\}.
    \end{equation}
    Therefore, if $\alpha \leq \frac{1}{\lambda_{\text{max}} G}$, we have:
    $$\alpha \leq \min _{i \in\left\{i \in I_k^-\mid g_{k}^{i}>0, p_{k}^{i}>0\right\}}
    \left\{\frac{x_{k}^{i}}{p_{k}^{i}}\right\}.$$
    From Lipschitz continuity of $\nabla f$, we have that for any $\alpha $ satisfying $\alpha \leq \frac{1}{\lambda_{\text{max}} G}$,
    \begin{equation}
        \begin{aligned}
            f\left(x_{k}(\alpha)\right)-f\left(x_{k}\right) & \leq g_{k}^{T}\left(x_{k}(\alpha)-x_{k}\right)+\frac{L}{2}\left\|x_{k}(\alpha)-x_{k}\right\|^{2} \\
            & \leq (\frac{\alpha^2 \lambda_{\text{max}}L}{2}-\alpha) \sum_{i \in I_k^-} g_k^i p_k^i + (\frac{\alpha \lambda_{\text{max}}L}{2}-1) \sum_{i \in I_{k}^{+}}  g_{k}^{i} (x_{k}^i-x_{k}^i(\alpha)).
        \end{aligned}
    \end{equation}
    If $\alpha$ satisfies $\alpha \leq \frac{2(1-\sigma)}{L \lambda_{\max }} $, then we have:
    \begin{equation}
        f\left(x_{k}\right)-f\left(x_{k}\left(\alpha\right)\right) \geq \sigma\left\{\alpha \sum_{i \notin I_{k}^{+}} g_k^i p_{k}^{i}
        +\sum_{i \in I_{k}^+} g_k^i (x_{k}^{i}-x_{k}^{i}\left(\alpha\right))\right\}.
    \end{equation}
    Let $\bar{\alpha} \triangleq \min \left\{\frac{1}{\lambda_{\max } G}, \frac{2(1-\sigma)}{L \lambda_{\max }}\right\}$.
    Due to the definition of $m_k$, $m_{k} \leq \min \left\{m \in \mathbb{N} \mid \beta^{m} \leq \bar{\alpha}\right\}$
    and $\beta^{m_{k}-1}>\bar{\alpha}$. Thus $m_{k}<+\infty$. This further indicates that
    \begin{equation}
        \alpha_{k} / \beta=\beta^{m_{k}-1}>\bar{\alpha} \Longrightarrow \alpha_{k}>\bar{\alpha} \beta.
    \end{equation}
    Next, we consider the iteration complexity of Algorithm 1.
    \begin{equation}
        \|S_k \nabla f(x_k)\|>\epsilon \Longleftrightarrow  \|S_k \nabla f(x_k)\|^2>\epsilon^2 \Longrightarrow \sum_{i \in I_k^+}(s_k^i g_k^i)^2 > \frac{\epsilon^2}{2} \text{ or } \sum_{i \in I_k^-}(s_k^i g_k^i)^2 > \frac{\epsilon^2}{2}
    \end{equation}
    If
    \begin{equation}
        \sum_{i \in I_k^+}(s_k^i g_k^i)^2 > \frac{\epsilon^2}{2} \Longleftrightarrow \sum_{i \in I_k^+}(x_k^i g_k^i)^2 > \frac{\epsilon^2}{2} \Longrightarrow
        \sum_{\substack{i \in I_k^+ \\ x_k^i - \alpha p_k^i \geq0}}(x_k^i g_k^i)^2 > \frac{\epsilon^2}{4} \text{ or }\sum_{\substack{i \in I_k^+ \\ x_k^i - \alpha p_k^i <0}}(x_k^i g_k^i)^2 > \frac{\epsilon^2}{4}
    \end{equation}
    Then,
    \begin{equation}
        \begin{aligned}
            \sum_{i \in I_{k}^+} g_k^i (x_{k}^{i}-x_{k}^{i}\left(\alpha\right))
            =&\sum_{\substack{i \in I_k^+ \\ x_k^i - \alpha p_k^i \geq0}} \alpha p_k^i g_k^i + \sum_{\substack{i \in I_k^+ \\ x_k^i - \alpha p_k^i <0}} x_k^i g_k^i\\
            = &\sum_{\substack{i \in I_k^+ \\ x_k^i - \alpha p_k^i \geq0}} \alpha (x_k^i g_k^i)^2 D_k[i,i] + \sum_{\substack{i \in I_k^+ \\ x_k^i - \alpha p_k^i <0}} x_k^i g_k^i\\
            \geq & \sum_{\substack{i \in I_k^+ \\ x_k^i - \alpha p_k^i \geq0}} \alpha \lambda_{\text{min}}(x_k^i g_k^i)^2 + \sqrt{\sum_{\substack{i \in I_k^+ \\ x_k^i - \alpha p_k^i <0}} (x_k^i g_k^i)^2}\\
            \geq & \min \left\{\bar{\alpha}\beta \lambda_{\text{min}}\frac{\epsilon^2}{4},\frac{\epsilon}{2}\right\}
            \stackrel{\epsilon<1}{\geq} \min \left\{\bar{\alpha}\beta \lambda_{\text{min}}\frac{\epsilon^2}{4},\frac{\epsilon^2}{2}\right\}
        \end{aligned}
    \end{equation}
    Else,
    \begin{equation}
        \sum_{i \in I_k^-}(s_k^i g_k^i)^2 > \frac{\epsilon^2}{2}
    \end{equation}
    Then,
    \begin{equation}
        \alpha \sum_{i \notin I_{k}^{+}} g_k^i p_{k}^{i}= \alpha \bar{g}_{k}^{T}\bar{p}_{k}=
        \alpha \bar{g}_{k}^{T}\bar{S}_k \bar{D}_{k} \bar{S}_k \bar{g}_{k} \geq \alpha \lambda_{\text{min}}\|\bar{S}_k \bar{g}_{k}\|^2
        \geq \bar{\alpha} \beta \lambda_{\text{min}}\|\bar{S}_k \bar{g}_{k}\|^2
        \geq \bar{\alpha} \beta \lambda_{\text{min}}\frac{\epsilon^2}{2}
    \end{equation}
    Therefore, if $\|S_k \nabla f(x_k)\|>\epsilon$ hold, from (22)(26)(28) we have:
    \begin{equation}
        \begin{aligned}
            f\left(x_{k}\right)-f\left(x_{k}\left(\alpha\right)\right) \geq & \sigma\left\{\alpha \sum_{i \notin I_{k}^{+}} g_k^i p_{k}^{i}
            +\sum_{i \in I_{k}^+} g_k^i (x_{k}^{i}-x_{k}^{i}\left(\alpha\right))\right\}\\
            \geq & \sigma \min \left\{\frac{\bar{\alpha}\beta \lambda_{\text{min}}}{4},\frac{1}{2}\right\}\epsilon^2
        \end{aligned}
    \end{equation}
    Suggest that the algorithm will stop within $K$ number of iterations, then,
    \begin{equation}
        f(x_0)-f_{\text{low}}\geq \sum_{k=0}^{K-1}f\left(x_{k}\right)-f\left(x_{k}\left(\alpha\right)\right)=f(x_0)-f(x_K)
        \geq K \sigma \min \left\{\frac{\bar{\alpha}\beta \lambda_{\text{min}}}{4},\frac{1}{2}\right\}\epsilon^2
    \end{equation}
    This further indicates that,
    \begin{equation}
        K \leq \frac{f(x_0)-f_{\text{low}}}{\sigma \min \left\{\frac{\bar{\alpha}\beta \lambda_{\text{min}}}{4},\frac{1}{2}\right\}\epsilon^2}
        = \frac{(f(x_0)-f_{\text{low}})\max \left\{4\lambda_{\max } G, \frac{2L \lambda_{\max }}{1-\sigma},2\beta \lambda_{\text{min}}\right\}}{\sigma \beta \lambda_{\text{min}} \epsilon^2}
    \end{equation}
\end{proof}

\section{Generalize the method for $\ell_1$-norm regularization problem}
We consider the $\ell_1$-norm regularization problem:
\begin{equation}
        \min \psi(x) \triangleq f(x)+\gamma \|x\|_1
\end{equation}
where  $f: \mathbb{R}^{n} \rightarrow \mathbb{R}$  is Lipschitz continuously differentiable 
and is bounded below by $f_{\text {low }}$. For any local optimal point of $\ell_1$-norm regularization problem, we give the following first-order necessary conditions:
\begin{theorem}[First-order necessary conditions]
    If $x^*$ is an local optimal point of $\ell_1$-norm regularization problem, then 
    \begin{equation}
        -\nabla f\left(x^{*}\right) \in \gamma \partial\left\|x^{*}\right\|_{1},
    \end{equation}
    specifically,
    \begin{equation}
    \label{1o-l1}
        \nabla_{i} f\left(x^{*}\right)=\left\{\begin{array}{ll}
-\gamma, & x_{i}^{*}>0, \\
a \in[-\gamma, \gamma], & x_{i}^{*}=0, \\
\gamma, & x_{i}^{*}<0 .
\end{array}\right.
    \end{equation}
\end{theorem}
%{\yx Are you able to show this?}

We propose an two-metric adaptive projection algorithm for resolution. We update $x_{k+1}$ using:
$$x_{k+1}:=\mathcal{P}_k\left(x_{k}-\alpha_{k} D_{k} (g_{k}+\omega_k)\right), \quad \forall k \geq 0,$$
where $\mathcal{P}_k$ and $\omega_k$ is associated to $x_k$ and defined as below:
\begin{equation}
    \mathcal{P}_k^i(y)= 
\begin{cases}
  & \max \left\{y^i,0\right\}\quad \text{ if } x_k^i>0 \text{ or } x_k^i=0, g_k^i \leq -\gamma \\
  & \min \left\{y^i,0\right\}\quad \text{ if } x_k^i<0 \text{ or } x_k^i=0, g_k^i \geq \gamma\\
  & 0 \quad \quad\quad\quad \quad\text{ if } x_k^i=0, |g_k^i|<\gamma  
\end{cases}
\end{equation}
\begin{equation}
    \omega_k= 
\begin{cases}
  & \gamma \quad \text{ if } x_k^i>0 \text{ or } x_k^i=0, g_k^i \leq -\gamma \\
  & -\gamma\quad \text{ if } x_k^i<0 \text{ or } x_k^i=0, g_k^i \geq \gamma\\
  & 0 \quad \text{ if } x_k^i=0, |g_k^i|<\gamma  
\end{cases}
\end{equation}
where $\alpha_{k}>0$ is the stepsize, 
$g_{k} \triangleq \nabla f\left(x_{k}\right)$. 
$D_{k}$ is a symmetric positive definite matrix in $\mathbb{R}^{n \times n}$.
Define $p_{k} \triangleq D_{k}  (g_{k}+\omega_k)$, $I_{k}^{+} \triangleq \left\{i\mid x_k^i=0, |g_k^i|<\gamma\right\}$, $I_{k}^{-} \triangleq[1, \ldots, n] \backslash I_{k}^{+}$.
We suppose that  $D_{k}$  is diagonal w.r.t.  $I_{k}^{+}$ , i.e.,
\begin{equation}
    D_{k}[i, j]=0, \quad \forall i \in I_{k}^{+}, j \neq i .
\end{equation}
Denote $x_{k}(\alpha) \triangleq P\left(x_{k}-\alpha p_{k}\right), \forall \alpha>0$.
Let  $\bar{g}_{k}$, $\bar{p}_{k}$ and $\bar{\omega}_k$ be the subvector of  $g_{k}$, $p_{k}$ and $\omega_k$  w.r.t index set  $I_{k}^{-}$ , and  $\bar{D}_{k}$ be the submatrix of $D_{k}$ and w.r.t.  $I_{k}^{-} \times I_{k}^{-}$. 
Then the algorithm can be formally stated as in Algorithm~\ref{alg: 2madapproj}.
\begin{algorithm}\label{alg: 2madapproj}
    \SetAlgoLined
    \KwIn{step acceptance parameter $\sigma \in(0,1)$, backtracking parameter $\beta \in(0,1)$, regularization parameter $\gamma$ and initial point $x_{0}$.}%输入参数
    \For{$k=0,\ldots,$ until achieve the stopping criteria}
    {
        Find:
            \begin{equation}
                m_k \triangleq \min \left\{m \in \mathbb{N} \Bigg| \psi \left(x_{k}\right)-\psi \left(x_{k}\left(\beta^{m}\right)\right) 
            \geq \sigma \beta^{m} (\bar{g}_k+\bar{\omega}_k) \bar{p}_{k}\right\}
            \end{equation}

        Let $\alpha_k= \beta^{m_k}$, $x_{k+1}= x_k(\alpha_k)$.
    }
    \KwOut{$x_k$}
    \caption{Two-metric adaptive projection for $\ell_1$-norm minimization}
\end{algorithm}

For Algorithm 2, we have the following three theorems.
\begin{theorem} 
$\forall \alpha>0,x_{k}(\alpha) = x_{k}$ if and only if $x_k$ satisfy (\ref{1o-l1}).
\end{theorem}
\begin{proof}
    Assume that for all $\alpha>0$, we have $x_{k}(\alpha) = x_{k}$. Then, if $x_k^i=0, g_k^i \geq \gamma$, we have:
    \begin{equation}
        x_k^i = \min \left\{x_k^i - \alpha p_k^i,0\right\}= \min \left\{x_k^i - \alpha [D_k(g_k+\omega_k)]^i,0\right\} = 0,
    \end{equation}
    which indicates that
    \begin{equation}
        x_k^i - \alpha p_k^i \geq 0 \Longleftrightarrow p_k^i \leq 0.
    \end{equation}
    Note that,
    \begin{equation}
        g_k^i+\omega_k^i = g_k^i-\gamma \geq 0,
    \end{equation}
    combine the above two equations, we get
    \begin{equation}
        (g_k^i+\omega_k^i)p_k^i \leq 0.
    \end{equation}
    For $x_k^i=0, g_k^i \leq -\gamma$, we have
    \begin{equation}
        p_k^i \geq 0 \text{ and } g_k^i+\omega_k^i \leq 0 \Longrightarrow (g_k^i+\omega_k^i)p_k^i \leq 0.
    \end{equation}
    For $x_k^i > 0$, we have 
    \begin{equation}
        x_k^i = \max \left\{x_k^i - \alpha p_k^i,0\right\} >0 \Longrightarrow x_k^i = x_k^i - \alpha p_k^i \Longrightarrow p_k^i=0.
    \end{equation}
    For $x_k^i < 0$, we also have
    \begin{equation}
        p_k^i=0.
    \end{equation}
    To sum up, we have 
    \begin{equation}
        \sum_{i \in I_k^-}(g_k^i+\omega_k^i)p_k^i \leq 0 \Longrightarrow (\bar{g}_k+\bar{\omega}_k) \bar{D}_k (\bar{g}_k+\bar{\omega}_k) \leq 0.
    \end{equation}
    Since $D_k$ is positive definite, we must have
    \begin{equation}
        \bar{g}_k+\bar{\omega}_k = 0
    \end{equation}
    which indicate (\ref{1o-l1}). The other side is trivial.
\end{proof}

\begin{theorem}
    For $\sigma \in (0,1)$, $\alpha$ small enough, we have
    \begin{equation}
        \psi \left(x_{k}\right)-\psi \left(x_{k}\left(\alpha\right)\right) 
            \geq \sigma \alpha (\bar{g}_k+\bar{\omega}_k) \bar{p}_{k}
    \end{equation}
\end{theorem}
\begin{proof}
    From Lipschitz continuity of $\nabla f$, we have the following inequality:
    \begin{equation}
        f(x_{k}(\alpha))-f(x_{k}) \leq g_{k}^{T} (x_{k}(\alpha)-x_{k})+\frac{L}{2}\|x_{k}(\alpha)-x_{k}\|_{2}^{2},
    \end{equation}
    so that
    \begin{equation}
        \psi(x_{k}(\alpha))-\psi(x_{k}) \leq g_{k}^{T} (x_{k}(\alpha)-x_{k})+\frac{L}{2}\|x_{k}(\alpha)-x_{k}\|_{2}^{2}+\gamma(\|x_{k}(\alpha)\|_1-\|x_{k}\|_1).
    \end{equation}
    For $x_k^i \ne 0$, $x_k^i$ and $x_k^i(\alpha)$ have the same sign for $\alpha$ small enough.
    So that we have
    \begin{equation}
        \begin{aligned}
           g_{k}^{T} (x_{k}(\alpha)-x_{k})+\gamma(\|x_{k}(\alpha)\|_1-\|x_{k}\|_1) & =\sum_{\substack{ x_k^i \ne 0}}(g_k^i+\sign{x_k^i}\gamma)(x_k^i(\alpha)-x_k^i)\\
           & + \sum_{\substack{ x_k^i = 0 \\ g_k^i \leq -\gamma}}(g_k^i+\gamma)x_k^i(\alpha)+\sum_{\substack{ x_k^i = 0 \\ g_k^i \geq \gamma}}(g_k^i-\gamma)x_k^i(\alpha)\\
            & = -\alpha \sum_{i \in I_k^-} (g_k^i+\omega_k^i)p_k^i+\sum_{\substack{i \in I_k^- \\ x_k^i = 0}} (g_k^i+\omega_k^i)(x_k^i(\alpha)+\alpha p_k^i)
        \end{aligned}
    \end{equation}
    Note that, for $x_k^i=0, g_k^i \leq -\gamma$, we have
    \begin{equation}
        g_k^i+\omega_k^i \leq 0 \text{ and } x_k^i(\alpha)+\alpha p_k^i\geq0,
    \end{equation}
    for $x_k^i=0, g_k^i \geq \gamma$, we have
    \begin{equation}
        g_k^i+\omega_k^i \geq 0 \text{ and } x_k^i(\alpha)+\alpha p_k^i \leq0.
    \end{equation}
    From the above two equations, we get 
    \begin{equation}
        \sum_{\substack{i \in I_k^- \\ x_k^i = 0}} (g_k^i+\omega_k^i)(x_k^i(\alpha)+\alpha p_k^i) \leq 0 \Longrightarrow g_{k}^{T} (x_{k}(\alpha)-x_{k})+\gamma(\|x_{k}(\alpha)\|_1-\|x_{k}\|_1) \leq -\alpha \sum_{i \in I_k^-} (g_k^i+\omega_k^i)p_k^i 
    \end{equation}
    For all  $\alpha \geq 0$, we also have
    \begin{equation}
        \left|x_{k}^{i}-x_{k}^{i}(\alpha)\right| \leq \alpha\left|p_{k}^{i}\right| \quad \forall i \in I_k^-,
    \end{equation}
    so that
    \begin{equation}
        \frac{L}{2}\|x_{k}(\alpha)-x_{k}\|_{2}^{2} \leq \frac{L \alpha^2}{2}\|\bar{p}_k\|_{2}^{2} = \frac{L \alpha^2}{2}\|\bar{D}_k(\bar{g}_k+\bar{\omega}_k)\|_{2}^{2} \leq \frac{L \alpha^2 \lambda_{\max}}{2}(\bar{g}_k+\bar{\omega}_k)\bar{p}_k.
    \end{equation}
    Overall, if $\alpha \leq \frac{2(1-\sigma)}{L \lambda_{\max}}$ and small enough such that $x_k^i$ and $x_k^i(\alpha)$ have the same sign when $x_k^i \ne 0$.,
    \begin{equation}
        \psi(x_{k}(\alpha))-\psi(x_{k}) \leq \left(\frac{L \alpha^2 \lambda_{\max}}{2}-\alpha\right)\left( \bar{g}_k+\bar{\omega}_k \right)\bar{p}_k \leq -\sigma \alpha (\bar{g}_k+\bar{\omega}_k)\bar{p}_k,
    \end{equation}
\end{proof}
\begin{theorem}
    The right-hand side of (50) is nonnegative and is positive for all $\alpha>0$ if and only if $x_{k}$ does not satisfy (\ref{1o-l1}).
\end{theorem}
\begin{proof}
    Note that
    \begin{equation}
        (\bar{g}_k+\bar{\omega}_k)\bar{p}_k = (\bar{g}_k+\bar{\omega}_k) \bar{D}_k (\bar{g}_k+\bar{\omega}_k) \geq 0,
    \end{equation}
    since $D_k$ is positive definite. Note that
    \begin{equation}
        x_{k} \text{ satisfy (\ref{1o-l1})} \Longleftrightarrow \bar{g}_k+\bar{\omega}_k = 0.
    \end{equation}
    Therefore, both sides are trivial.
\end{proof}
In conclusion, the algorithm is well-defined and decreases the value of the objective function at each iteration $k$, in which $x_{k}$ is not a critical point.

\subsection{Numerical experiment}
We use Algorithm 2 to solve the LASSO problem:
$
    \min \left\{ \frac{1}{2}\|Ax-b\|_2^2+\gamma \|x\|_1\right\},
$
where $A \in \mathbb{R}^{m\times n}$$(m \ll n)$ and each element of $A$ follows the normal distribution, 
and the real solution $u \in \mathbb{R}^n$ is a sparse vector with only 10\% none zero element, we derive $b$ by setting $b = Au$.
To accelerate the algorithm's convergence rate, a continuity strategy can be adopted to gradually reduce the larger regularization parameter from $\gamma_0$ to $\gamma$.
Then the algorithm can be formally stated as follows.

\begin{algorithm}[ht]
    \SetAlgoLined
    \KwIn{step acceptance parameter $\sigma \in(0,1)$, reduction factor $\eta \in(0,1)$, regularization parameter $\gamma_0, \gamma$ and initial point $x_{0}$. Let $k=0$.}%输入参数
    \While{achieving the stopping criteria}
    {
        $x_{k+1} \longleftarrow$ solution of Algorithm 2 for LASSO with $\gamma_k$ and $x_k$.

        $\gamma_{k+1} \longleftarrow \max \left\{\eta \gamma_{k+1}, \gamma \right\}$

        $k \longleftarrow k+1$
    }
    \KwOut{$x_k$}
    \caption{Two-metric adaptive projection for LASSO}
\end{algorithm}
Note that, $f(x)$ in LASSO is $f(x)=\frac{1}{2}\|Ax-b\|_2^2$ and $\nabla^2 f(x)=A^{T}A$, we set $D_k=A^{T}A + \epsilon I$ in algorithm 2, since $\nabla^2 f(x)$ is positive semi-definite.
The result is shown in Figure 1.
\begin{figure}
    \centering
    \includegraphics[scale=0.3]{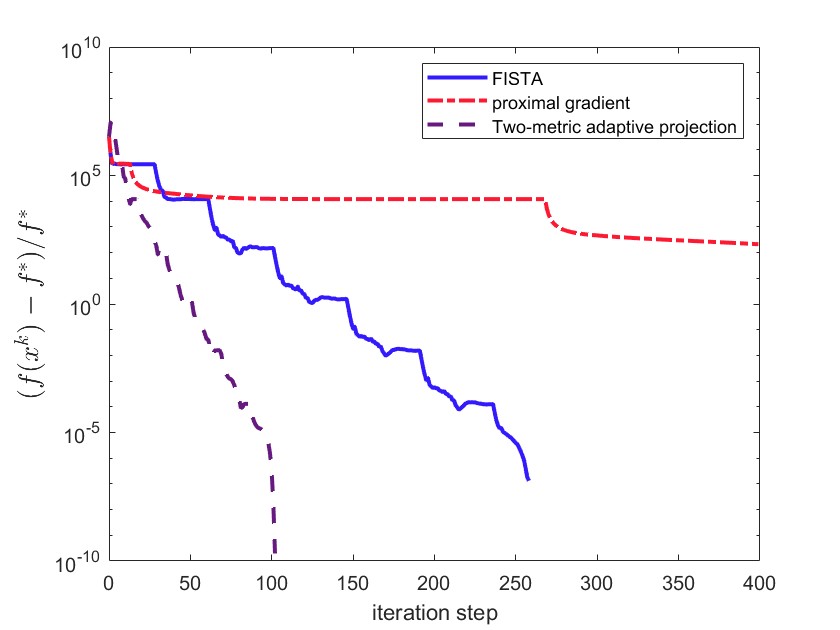}
    \caption{Comparison of methods on LASSO}
\end{figure}
The result shows that two-metric adaptive projection may have local superlinear convergence rate if we use Hessian, in sharp contrast with first-order methods such as FISTA \cite{doi:10.1137/080716542} and proximal gradient \cite{parikh2014proximal}.
\section{Conclusion}
In this article, we first define an approximate first-order optimal point for the bound-constrained problem. We investigate the iteration complexity of the two-metric projection method proposed by Bertsekas, the algorithm terminates within $\mathcal{O}(\epsilon^{-3})$ iterations and finds a point
that is approximately first-order optimal to tolerance $\epsilon$. By scaling the diagonal matrix, we equip the method
with a worst-case complexity theory that matches the best-known theoretical bounds for bound-constrained optimization and even unconstrained optimization. We also try to generalize the two-metric projection algorithm for composite optimization, we propose a framework called two-metric adaptive projection, which is suitable for solving the $\ell_1$-norm minimization problem. 

In future work, we will consider the convergence properties of our two-metric adaptive projection and propose a general algorithm for composite optimization.
\bibliographystyle{plain}
\bibliography{reference}

\end{document}